\input amstex
\documentstyle{amsppt}
\input epsf
\magnification1100
\TagsOnLeft
\parskip0pt plus1pt
\normalbaselineskip12pt
\normallineskip2pt
\normallineskiplimit1pt
\normalbaselines
\pagewidth  {6.50 truein}
\pageheight {9.00 truein}
\hcorrection{0.00 truein}
\vcorrection{0.00 truein}
\NoPageNumbers
\BlackBoxes
\def\parindentamount{20pt}
\parindent\parindentamount
\nologo
\def\CmD{\C \setminus \overline D}
\def\dotfil{\leaders\hbox to0.5em{\hss.\hss}\hfil}
\def\hfilll{\hskip 0pt plus 1filll}

\def\C{{\Bbb C}}
\def\R{{\Bbb R}}
\def\Z{{\Bbb Z}}
\def\cf{{\it cf.}}
\def\ie{{\it i.e.}}
\def\set#1#2{\left\{\,#1\,\left\vert\,#2\right.\,\right\}} 
\def\hangpar{\hangindent=0.45truein\hangafter=1}

\document
\centerline{COMPLEX HORSESHOES AND THE DYNAMICS}
\centerline{OF MAPPINGS OF TWO COMPLEX VARIABLES}
\bigskip
\centerline{A Thesis}
\centerline{Presented to the Faculty of the Graduate School}
\centerline{of Cornell University}
\centerline{in Partial Fulfillment of the Requirements for the Degree of}
\centerline{Doctor of Philosophy}
\bigskip
\centerline{by}
\centerline{Ralph Werner Oberste-Vorth}
\centerline{August, 1987}

\bigskip\hrule\bigskip

\centerline{\copyright\ Ralph Werner Oberste-Vorth 1987}
\centerline{ALL RIGHTS RESERVED}

\bigskip\hrule\bigskip

\centerline{COMPLEX HORSESHOES AND THE DYNAMICS}
\centerline{OF MAPPINGS OF TWO COMPLEX VARIABLES}
\bigskip
\centerline{Ralph Werner Oberste-Vorth, Ph.D.}
\centerline{Cornell University 1987}
\bigskip

In this study, a theory analogous to both the theories of polynomial-like mappings and Smale's real horseshoes is developed for the study of the dynamics of mappings of two complex variables.

In partial analogy with polynomials in a single variable there are the H\'enon mappings in two variables as well as higher dimensional analogues. From polynomial-like mappings, H\'enon-like and
quasi-H\'enon-like mappings are defined following this analogy. A special form of the latter is the complex horseshoe.

The major results about the real horseshoes of Smale remain true in the complex setting. In particular:

\bigskip\hrule\bigskip

\roster
 \item Trapping fields of cones(which are sectors in the real case) in the tangent spaces can be defined and used to find horseshoes.
 \item The dynamics of a horseshoe is that of the two-sided shift on the symbol space on some number of symbols which depends on the type of the horseshoe.
 \item Transverse intersections of the stable and unstable manifolds of a hyperbolic periodic point guarantee the existence of horseshoes.
 \endroster

\bigskip\hrule\bigskip

\centerline{Biographical Sketch}
\bigskip

Ralph Oberste-Vorth was born on May 13, 1959 in Brooklyn, New York. He graduated from Stuyvesant High School in New York in 1977 and then attended Hunter College and the Graduate School and University
Center of the City University of New York, majoring in mathematics with a minor in computer science. In 1981, he received the Bachelor of Arts (Magna Cum Laude) and Master of Arts degrees from Hunter
College. Since 1981 he has been a graduate student in mathematics at Cornell University. He has held adjunct faculty positions at Hunter College and John Jay College of Criminal Justice of CUNY and the
State University of New York College at Cortland.

\bigskip\centerline{iii}\bigskip\hrule\bigskip

\centerline{To my parents}

\bigskip\centerline{iv}\bigskip\hrule\bigskip

\centerline{Acknowledgments}
\bigskip

I wish to thank Professor John Hamal Hubbard for his guidance, inspiration, confidence, and friendship. At a time when I thought mathematics had died, his work made it come back to life for me. His energy
and accessibility, day and night, have been very much appreciated.

I thank Professors Adrien Douady, John Guckenheimer, and John Smillie for their interest and comments. I thank Professor Clifford Earle for serving on my Special Committee for the entire duration of my stay
at Cornell.

I thank my parents, to whom this is dedicated, for their patience and encouragement despite their doubts.

I thank Homer Smith for his expert computer and photographic assistance and for generously loaning me some of his equipment.

Last, but by no means least, I thank my friends for providing the many diversions, including bridge, squash, and tennis, as well as the comic relief needed to make life bearable.

\bigskip\centerline{v}\bigskip\hrule\bigskip

\centerline{Table of Contents}
 $$
 \vbox{\halign{\hfil#&\quad\hbox to 4truein{#\dotfil}&\ \  Page\
 \hfil#\cr 
    &Introduction &1\cr 
 1. &H\'enon-like mappings &8\cr 
 2. &Complex horseshoes &20\cr 
 3. &Complex horseshoes as shift dynamical systems &37\cr 
 4. &Ubiquity of complex horseshoes &46\cr 
    &List of References &54\cr
 }}
 $$

\bigskip\centerline{vi}\bigskip\hrule\bigskip

 \centerline{Table of Tables}
 $$
 \vbox{\halign{\hfil#&\quad\hbox to 4truein{#\dotfil}&\ \  Page\
 \hfil#\cr 
 1.1 &Periodic cycles of H\'enon-like mappings of degree 2 &17\cr
 }}
 $$

\bigskip\centerline{vii}\bigskip\hrule\bigskip

 \centerline{Table of Figures}
 $$
 \vbox{\halign{\hfil#&\quad\hbox to 4truein{#\dotfil}&\ \  Page\
 \hfil#\cr 
 1.1 &H\'enon-like mappings of degree 2 &9\cr
 4.1 &Horseshoes from transverse homoclinic points &50\cr
 }}
 $$

\bigskip\centerline{viii}\bigskip\hrule\bigskip

\centerline{\bf Introduction}\bigskip

The study of the subject of the dynamics of complex analytic functions of one variable goes back to the early 1900's with the publication of several papers by Pierre Fatou and Gaston Julia,
including the long memoirs \cite{F1} and \cite{J}, published around 1920. These papers proved to be the definitive works in the subject for a long time to come. For the most part, new interesting results
did not appear again until 1982. (For a survey of the subject, see \cite{Bl}.)

Somewhat surprisingly, this recent resurgence of the field has been spurred not so much by mathematical developments but by advances in computer graphics. However the original impetus at the beginning of
the century was strictly mathematical. Besides the obvious connection between differential equations and dynamical systems, the theory of complex analytic dynamics was very much based upon the developments
in complex analysis. In particular, the theory of normal families and the work of P. Montel in that regard had given rise to much of the work of Fatou and Julia.

Compared with the theory of complex analytic functions of a 

\bigskip\centerline{1}\bigskip\hrule\bigskip

\noindent
single variable, the theory of complex analytic mappings of several, that is, two or more, variables is quite different. In particular, the theory
of omitted values, including normal families, is not paralleled in the several complex variables case. This difference was really exposed by Fatou and L. Bieberbach in the 1920's in \cite{F2} and \cite{Bi}.

They showed the existence of what we refer to as {\it Fatou-Bieberbach domains}: open subsets of $\C^n$ whose complements have nonempty interior and yet are the images of $\C^n$ under an injective analytic
mapping. This is contrary to the one variable case where the image of every non-constant analytic function on $\C$ omits at most a single point.

The present work started with an attempt to understand such Fatou-Bieberbach domains. These arise naturally as the basins of attractive fixed points of analytic automorphisms of $\C^n$. The basins are then
the image of the mapping conjugating the given automorphism to its linear part at the given fixed point. 

This remains true even when the Jordan canonical form of the linearization is not diagonal and when there are resonances in general. This latter result is apparently new. (See \cite{HO}.) For 

\bigskip\centerline{2}\bigskip\hrule\bigskip

\noindent 
example, consider
 $$
 F: \bmatrix x \\ y \endbmatrix \mapsto \bmatrix x^2 + 9/32 - y/8 \\ x \endbmatrix.
 $$ 
 This has two fixed points, of which $(3/8, 3/8)$ is attractive with its linear part having resonant eigenvalues $1/4$ and $1/2$ (that is, $1/4 = (1/2)^2)$. Moreover, none of the points in the region
 $$
 \set{(x,y)}{|y| < 4 |x|^2/3, |x| > 4}
 $$
 remain bounded under iteration of $F$. So the basin of $(3/8,3/8)$ is not all of $\C^2$.

As was the case with one dimensional complex analytic dynamic, two dimensional complex dynamics had been a rather dormant field until recently. In analogy, the present work started with a computer
investigation of a specific Fatou-Bieberbach domain.

In the time between Fatou and the present, most of the attention of those studying dynamical systems has been limited to mappings in the real. This is somewhat surprising for two reasons. First, small
perturbations of the coefficients of polynomial terms of 

\bigskip\centerline{3}\bigskip\hrule\bigskip

\noindent 
real mappings are liable to have large effects. For example, the number and periods of the periodic cycles may change. In the complex, the behavior
is more uniform. Second, the major tools of complex analysis do not apply. These include the theory of normal families and the naturally contracting Poincar\'e metric together with the contracting mapping
fixed point theorem. (Of course, in the real, there is always the restriction to maps with negative schwarzian derivative.)

Over the past 15 years, the two (real) parameter family of mappings, $F_{a,c}: \R^2 \to \R^2$, of the type 
 $$
 F_{a,c}: \bmatrix x \\ y \endbmatrix \mapsto \bmatrix x^2 + c - ay \\ x \endbmatrix,
 $$ 
 with $a \ne 0$, has received much attention. H\'enon, in \cite{H1} and \cite{H2}, first studied these numerically and they have become known as the {\it H\'enon mappings}. This family contains, up to
conjugation, most of the most interesting of the simplest nonlinear polynomial mappings of two variables. However, H\'enon mappings are still rather poorly understood and indeed the original question
concerning the existence of a strange attractor for any values of the parameters is still 

\bigskip\centerline{4}\bigskip\hrule\bigskip

\noindent 
unresolved today.

Despite the differences between the real and complex theories and the one variable and several variable theories, much of the development of the subject of complex analytic dynamics in several variables has
been conceived through analogy.

Recently, H\'enon mappings have started to be examined in the complex, that is, with both the variables and the parameters being complex (\cite{HO}, \cite{FM}). Besides this exposition, a joint paper with
John Hubbard \cite{HO} detailing the part of this work which concerns H\'enon mappings is in preparation. It has also been noted that there exist analogous mappings of higher degree of the form 
 $$
 \bmatrix x \\ y \endbmatrix \mapsto \bmatrix p(x) - ay \\ x \endbmatrix,
 $$ 
 where $p$ is a polynomial of degree at least two and $a \ne 0$ \cite{FM}. Note that these are always invertible with inverses given by 
 $$
 \bmatrix x \\ y \endbmatrix \mapsto \bmatrix y \\ (p(y) - x)/a \endbmatrix.
 $$

\bigskip\centerline{5}\bigskip\hrule\bigskip

\noindent 
Since H\'enon mappings are of this form, these are called the generalized H\'enon mappings. For polynomials $p$ of degree $d$, they are called H\'enon mappings of degree $d$.

The major development in the theory of real two dimensional dynamical systems has been the {\it horseshoe mapping\/} of Smale \cite{S}. The main contribution of this thesis is the study of the generalization
of horseshoe mappings to the complex setting.

Recently, A. Douady and J. Hubbard \cite{DH} developed the theory of {\it polynomial-like mappings}. Pol\-y\-nom\-i\-al-like mappings are designed to capture, in their topology, the essence of polynomial
mappings. These have helped to inspire an analogous class of mappings of several variables. 

The following is a summary of the contents of this thesis:

Chapter 1 gives the definition of {\it H\'enon-like mappings of degree $d$}, which are intended to capture the dynamical essence of H\'enon mappings, in analogy with polynomial-like mappings. The
relationship between H\'enon, H\'enon-like, and polynomial-like mappings is explored.

\bigskip\centerline{6}\bigskip\hrule\bigskip

Chapter 2 gives the definition of {\it complex horseshoes}, the complex analog of Smale's real horseshoes, as a special class of {\it quasi-H\'enon-like mappings}, which themselves generalize H\'enon-like
mappings. A criterion called a {\it trapping field of cones\/} is given which guarantees horseshoes. This is used to show that a large number of H\'enon mappings are complex horseshoes.

Chapter 3 shows that, in analogy with the real case, complex horseshoes of degree $d$ are conjugate to the two-sided shift on symbol space on $d$ symbols. 

Chapter 4 shows that complex horseshoes are ubiquitous in the subject of two dimensional complex analytic dynamics by showing that they appear almost every time a mapping has a hyperbolic periodic point.

\bigskip\centerline{7}\bigskip\hrule\bigskip

\centerline{\bf 1. H\'enon-Like Mappings}\bigskip 

In this chapter, a class of mappings is defined which attempt to capture in their topology the dynamics of H\'enon mappings of all degrees. These mappings will be called {\it H\'enon-like}. The definition
of H\'enon-like mappings is inspired by the definition of polynomial-like mappings (see \cite{DH}), which was designed to capture the topological essence of polynomials on some disc or, more generally, on
some open subset of $\C$ isomorphic to a disc.

A polynomial-like mapping of degree $d$ is a triple $(U, U', f)$, where $U$ and $U'$ are open subsets of $\C$ isomorphic to discs, with $U'$ relatively compact in $U$, and $f: U' \to U$ analytic and
proper of degree $d$. Note that it is convenient to think of polynomial-like of degree $d$ as meaning an analytic mapping $f: U \to \C$ such that $f(\partial U) \subset \C \setminus \overline U$ and
$f|_{\partial U}$ of degree $d$ (\cf\ \cite{DH}).

Figure 1.1 gives examples of the behavior, pictured in $\R^2$, which should be captured by the definition of H\'enon-like mappings of degree $2$. (In each case, the crescent-shaped region is the image of
the square with $A'$ the image of $A$, etc.)

\bigskip\centerline{8}\bigskip\hrule\bigskip

It seems clear that the behaviors described by (a) and (b) versus (c) and (d) in Figure 1.1 must be described differently, albeit analogously, trading ``horizontal'' for ``vertical.''
 
We set some notation: $d$ will always be an arbitrary fixed integer greater than one. Let $\pi_1, \pi_2: \C^2 \to \C$ be the projections onto the first and second coordinates, respectively. We will consider
a bidisc $B = D_1 \times D_2 \subset \C^2$, where $D_1, D_2 \subset \C$ are discs.  Vertical and horizontal ``slices" of $B$ are denoted by

\bigskip
 \centerline{\epsfbox{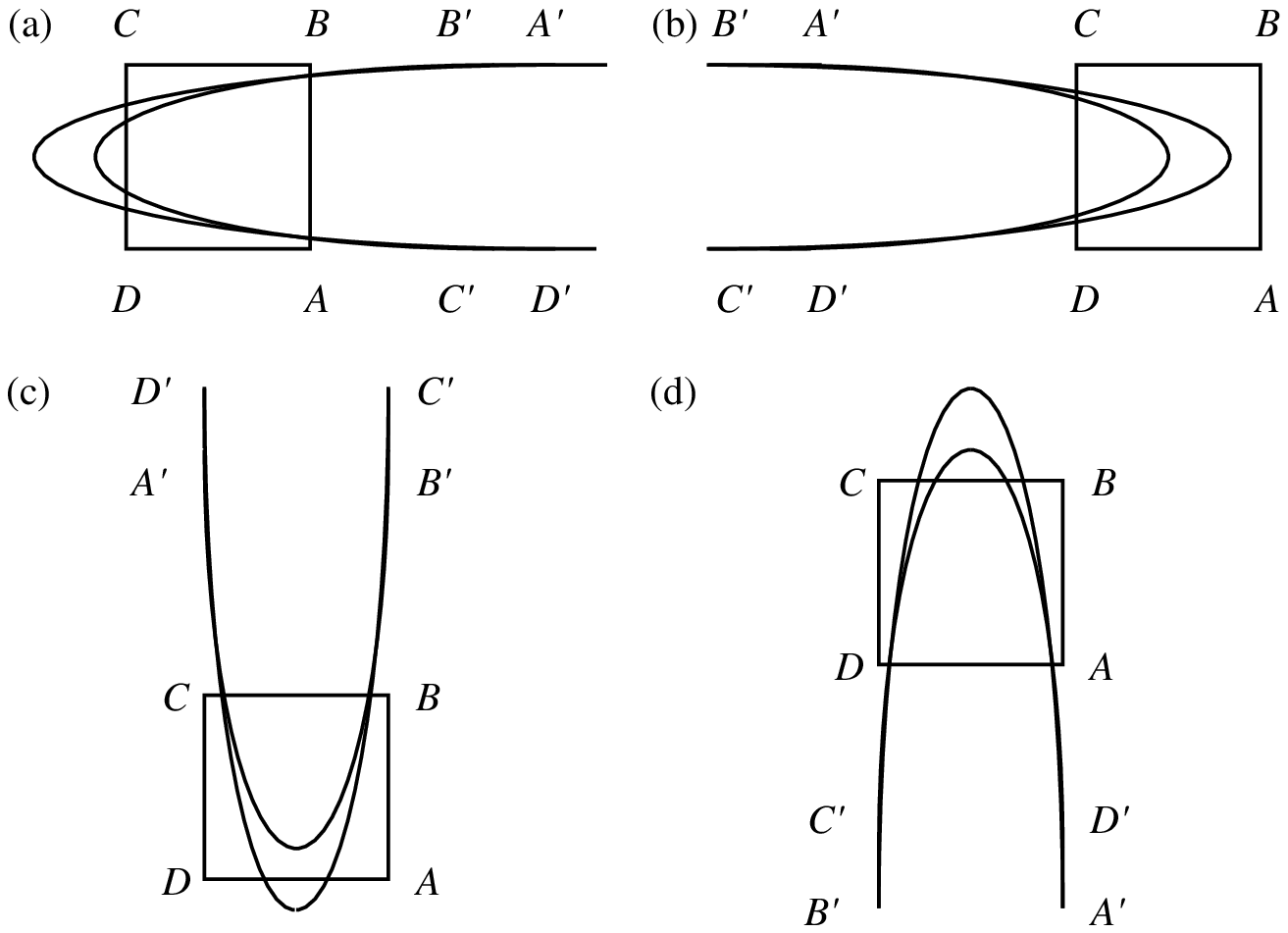}}
 \smallskip
 \centerline{Figure 1.1:{\quad}H\'enon-like of degree 2}
 \medskip

\bigskip\centerline{9}\bigskip\hrule\bigskip

$$
 V_x = \{x\} \times D_2 \qquad \text{and} \qquad H_y = D_1 \times \{y\},
 $$ 
 for all $x \in D_1$ and for all $y \in D_2$, respectively. We will be considering mappings of the bidisc, $F: \overline B \to \C^2$ together with a mapping denoted by $F^{-1}: \overline B \to \C^2$ which
is the inverse of $F$ where that makes sense. Now define, for each $(x, y) \in B$,  
 $$
 \align  
 F_{1,y} &= \pi_1 \circ F \circ (\text{Id} \times y): D_1 \to \C,\\  
 F_{2,x}^{-1} &= \pi_2 \circ F^{-1} \circ (x \times \text{Id}): D_2 \to \C,\\
 F_{2,x} &= \pi_2 \circ F \circ (x \times \text{Id}): D_2 \to \C,\\ 
 F_{1,y}^{-1} &= \pi_1 \circ F^{-1} \circ (\text{Id} \times y): D_1 \to \C. \qquad\hphantom{\text{and}} 
 \endalign
 $$

\definition {Definition 1.1} 
 $F: \overline B \to \C^2$ is a \underbar{\it H\'enon-like mapping of degree $d$} if there exists a mapping $G: \overline B \to \C^2$ such that 
 \roster 
 \item"(1)" Both $F$ and $G$ are injective and continuous on $\overline B$ and analytic on $B$. 
 \item"(2)" $F \circ G = \/{\roman Id}$ and $G \circ F = \/{\roman Id}$ where each makes sense.\newline
 Hence, rename $G$ as $F^{-1}$.  
 \item"(3)" For all $x \in D_1$ and $y \in D_2$, either 
 \item"" {\rm (a)}\quad $F_{1,y}$ and $F_{2,x}^{-1}$ are polynomial-like of degree $d$, or 
 \item"" {\rm (b)}\quad $F_{2,x}$ and $F_{1,y}^{-1}$ are polynomial-like of degree $d$. 
 \endroster 
 Depending on whether $F$ satisfies condition (a) or (b), call it

\bigskip\centerline{10}\bigskip\hrule\bigskip

\noindent 
\underbar{\it horizontal\/} or \underbar{\it vertical}.
 \enddefinition

\remark {Remarks} 
 (1) It would be correct and, perhaps, cleaner to define a H\'enon-like {\it pair\/} as a triple $(F, G, B)$ and then call $F$ a H\'enon-like mapping. However, this would put the emphasis in the wrong place.

(2) When the degree of a H\'enon-like mapping is either clear from context, or is not of primary interest, or is two, we will sometimes fail to mention the degree.

(3) Note that conditions (a) and (b) of (3) of Definition 1.1 are dual in the sense that if a H\'enon-like mapping $F$ satisfies condition (a), then $F^{-1}$ satisfies condition (b) and vice versa. In
Figure 1.1, (a) and (b) correspond with horizontal H\'enon-like mappings while (c) and (d) correspond with vertical H\'enon-like mappings. In general, unless otherwise specified, it will be assumed
that H\'enon-like mappings are horizontal, \ie, satisfying condition (a)

(4) Just as for polynomial-like mappings, it does not make sense, in general, to say that a mapping $F: \C^2 \to \C^2$ is H\'enon-like as this may be ambiguous. $F$ may exhibit different H\'enon-like
behavior in different regions. The domain $\overline B$ is part of the 

\bigskip\centerline{11}\bigskip\hrule\bigskip

\noindent 
definition of H\'enon-like mappings. Of course, the context may be used to resolve such ambiguities.

(5) Let
 $$
 \partial B_V = \partial \overline D_1 \times \overline D_2 \qquad \text{and} \qquad \partial B_H = \overline D_1 \times \partial \overline D_2.
 $$
 the ``vertical and horizontal boundaries."
 \endremark

\proclaim {Proposition 1.2}  If $F: \overline B \to \C^2$ is a H\'enon-like mapping, then either 
 $$
 \align 
 F (\partial B_V) \subset \C^2 \setminus \overline B \qquad &and \qquad F^{-1}(\partial B_H) \subset \C^2 \setminus \overline B \\ 
 \intertext{or}
 F^{-1}(\partial B_V) \subset \C^2 \setminus \overline B \qquad &and \qquad F(\partial B_H) \subset \C^2 \setminus \overline B.
 \endalign
 $$ 
 \endproclaim

\demo {Proof} 
 This follows from the fact that the boundary of a polynomial-like mapping is mapped outside of the closure of its domain.
 \newline\null\hfilll{\bf QED}
 \enddemo

Note that these are equivalent to

\bigskip\centerline{12}\bigskip\hrule\bigskip 

$$
 \align 
 F(\partial B_V) \cap \overline B = \emptyset \qquad&\text{and}\qquad F^{-1}(\partial B_H) \cap \overline B = \emptyset\\ 
 \intertext{or} 
 F^{-1}(\partial B_V) \cap \overline B = \emptyset \qquad&\text{and}\qquad F(\partial B_H) \cap \overline B = \emptyset.
 \endalign
 $$ 
 In the next chapter we will consider an alternative to H\'enon-like mappings where conditions (a) and (b) of Definition 1.1(3) are replaced by such conditions on the vertical and horizontal boundaries.

(6) The class of polynomial-like mappings is stable under small perturbations \cite{DH} and the same is true for H\'enon-like mappings:

 \proclaim {Proposition 1.3} 
 Suppose $F: \overline B \to \C^2$ is H\'enon-like of degree $d$. Let $H: \overline B \to \C^2$ be injective, continuous on $\overline B$, and analytic on $B$. If\/ $\|H\|$ is sufficiently small, then $F +
H$ is also H\'enon-like of degree $d$.
 \endproclaim

\demo {Proof} 
 Recall that $\partial B_V = \partial \overline D_1 \times \overline D_2$. By Proposition 1.2, we can choose an $\varepsilon > 0$ so that $\varepsilon < d(\overline B, F(\partial B_V))$. If $\|H(x, y)\| <
\varepsilon$ for all $(x, y) \in \partial B_V$, then $F + H$ is H\'enon-like of degree $d$.
 \newline\null\hfilll{\bf QED}
 \enddemo

\bigskip\centerline{13}\bigskip\hrule\bigskip

\noindent 
Of course, the H\'enon mappings themselves provide the obvious examples of H\'enon-like mappings:

\proclaim {Proposition 1.4} 
 For every H\'enon mapping, $F$, of every degree, there exists an $R$, such that $F: \overline D_R{}^2 \to \C^2$ is a H\'enon-like mapping of the same degree. 
 \endproclaim

\demo {Proof} 
 This follows immediately from the fact that all polynomials are polynomial-like of the same degree on sufficiently large discs. 
 \newline\null\hfilll{\bf QED}
 \enddemo

\example {Examples} 
 (1) A simple computation shows that if $F$ is the H\'enon mapping (of degree $2$) with parameters $a$ and $c$, and $D_R$ is the disc of radius $R$, with
 $$
 R > (1/2) \bigl(1 + |a| + \sqrt{(1 + |a|)^2 + 4|c|}\bigr),
 $$
 then $F: \overline D_R^2 \to \C^2$ is a horizontal H\'enon-like mapping of degree $2$. This $R$ is exactly what is required so that $F(\partial \overline D_R \times \overline D_R) \cap \overline D_R{}^2 =
\emptyset$ and $F^{-1}(\overline{D_R} \times \partial \overline{D_R}) \cap \overline{D_R}^2 = \emptyset$. Of course, their inverses are vertical H\'enon-like mappings.

\bigskip\centerline{14}\bigskip\hrule\bigskip

\noindent 
(2) More generally, consider the mappings $G: \C^2 \to \C^2$ of the form
 $$
 G: \bmatrix x \\ y \endbmatrix \mapsto \bmatrix x^d + c - ay \\ x \endbmatrix,
 $$
 with $a \ne 0$ and $d \ge 2$. When $d = 2$, we are back to the previous example. The lower bound on $R$ came from solving the inequality 
 $$
 R^d - (1 + |a|)R - 4|c| > 0
 $$ 
 for $R$. Note that when $d = 2$ we already had $R > 1$. Therefore, the same lower bound will work here as well. Of course, better lower bounds can be found.
 \endexample

Analogous to the invariant sets defined for H\'enon mappings \cite{HO}, we can define the following sets for H\'enon-like mappings:
 $$
 \gather
 {\align 
 K_+ &= \bigl\{\,z \in B\bigm| ^{\hphantom{-}}F^{\circ n}(z) \in B\quad \text{for all}\ n > 0\,\bigr\},\\ 
 K_- &= \bigl\{\,z \in B\bigm| F^{\circ -n}(z) \in B\quad \text{for all}\ n > 0\,\bigr\},
 \endalign}\\ 
 J_\pm = \partial K_\pm, \quad \ K = K_+ \cap K_-, \quad J = J_+ \cap J_-.
 \endgather
 $$

At least some of the dynamics of H\'enon mappings is captured by the notion of a H\'enon-like mapping. When the H\'enon-like mapping

\bigskip\centerline{15}\bigskip\hrule\bigskip

\noindent 
is actually a H\'enon mapping, as in Proposition 1.4, a simple computation
shows that the set $K$ for the H\'enon mapping actually is contained in $\overline{D_R}^2$ and thus equals the corresponding $K$ defined above. Moreover:

\proclaim {Proposition 1.5} 
 For every $d$, all H\'enon-like mappings of degree $d$ have the same number of periodic cycles, counted with multiplicity, as a polynomial of degree $d$. 
 \endproclaim

In the following, in order to simplify the language, the phrase ``counted with multiplicity'' is omitted but implied and the word ``period'' means the lowest period.

\demo {Proof} 
 Let $F: \overline B \to \C^2$ be an arbitrary H\'enon-like mapping of degree $d$, with coordinate functions $f_1 = \pi_1 \circ F$ and $f_2 = \pi_2 \circ F$. We may suppose that $F$ is horizontal. Consider
the family of mappings, $F_\varepsilon: \overline B \to \C^2$, for $0 < \varepsilon < 1$, defined by 
 $$
 F_\varepsilon(x, y) = (f_1(x, y) , \varepsilon f_2(x, y)).
 $$ 
 Each of the mappings $F_\varepsilon$ has the same number of periodic cycles. As $\varepsilon \to 0$, the periodic points converge to the plane $\C \times \{0\}$. 

\bigskip\centerline{16}\bigskip\hrule\bigskip

\noindent 
Hence we can consider the map on this plane induced by the limit mapping, $x \mapsto f_1(x, 0)$. By definition, this mapping is polynomial-like of degree $d$ on $H_0 = \{\,x\mid (x, 0) \in B\}$. Such a
mapping is quasi-conformally equivalent to a polynomial of degree $d$ and all polynomials of degree $d$ have the same number of periodic cycles. 
 \newline\null\hfilll{\bf QED}
 \enddemo

\remark {Remarks} 
 Moreover, all H\'enon mappings of degree $d$ have the same number of periodic cycles. In particular, they can be counted explicitly.

There are $d$ fixed points. The number of periodic cycles of prime period $p$ is $d(d^{p-1} - 1)/p$. For period $n \ge 4$, with $n$ not prime, there is a recursive algorithm for determining the number of
periodic points: the sum, taken over all $m$ which divide $n$, of the number of points of period $m$ equals $d^n$. Table 1.1 shows the 

\smallskip \goodbreak\midinsert
 \parindent0pt\hfil\vbox{\tabskip=0pt\offinterlineskip
 \def\tablerule{\noalign{\hrule}}
 \halign to 4.5in{\strut#& \vrule#\tabskip=1em plus2em&
  \hfil#\hfil& \vrule#& \hfil#\hfil& \vrule#&
  \hfil#\hfil& \vrule#& \hfil#\hfil& \vrule#&
  \hfil#\hfil& \vrule#& \hfil#\hfil& \vrule#\tabskip=0pt\cr\tablerule
 &&\omit Period&&\omit Cycles&&
  \omit Period&&\omit Cycles&&
  \omit Period&&\omit Cycles&\cr\tablerule
 &&1&&2&&5&& 6&& 9&& 56&\cr\tablerule
 &&2&&1&&6&& 9&&10&& 99&\cr\tablerule
 &&3&&2&&7&&18&&11&&186&\cr\tablerule
 &&4&&3&&8&&30&&12&&335&\cr\tablerule
 }}\hfil
 \smallskip\nobreak
 \centerline{Table 1.1:\ \ Periodic cycles of H\'enon-like mappings of degree two}
 \parindent\parindentamount\endinsert 

\bigskip\centerline{17}\bigskip\hrule\bigskip

\noindent 
number of periodic cycles of H\'enon-like mappings of degree $2$ for periods $1$ through $12$.
 \endremark

Of course, there are other examples of H\'enon-like mappings besides actual H\'enon mappings and their perturbations. This is much the same as for polynomial-like mappings.

It is tempting to conjecture that every H\'enon-like mapping of degree $2$, or, in fact, any degree, is conjugate to a H\'enon mapping, as in the case of polynomial-like mappings. Unfortunately, this is
false.

\example {Example} 
 Find an actual H\'enon mapping $F$ with $a = 1$ and two periodic points, say of period $k$, at each of which the linearization of $F^{\circ k}$ has complex conjugate eigenvalues, of absolute value one.
Then choose $R > 1 + \sqrt{1 + |c|}$ and a mapping
 $$
 H: \overline{D_R}^2 \to \C^2
 $$ 
 which is small on $\overline{D_R}^2$, vanishes at the two periodic points, and so that 
 $$
 F + H: \overline{D_R}^2 \to C^2
 $$ 

\bigskip\centerline{18}\bigskip\hrule\bigskip

\noindent 
is a H\'enon-like mapping with these points still periodic, but one being attractive and the other being repulsive. This cannot be conjugate to a H\'enon mapping.
 \endexample

However, we cannot find such a H\'enon mapping $F$ where the periodic points are fixed points are fixed points even if $a = 1$ is replaced with an arbitrary complex parameter $a$ with $|a| = 1$. The problem
is that when a H\'enon mapping has distinct fixed points, then at least one of them is hyperbolic.

Perhaps with the added requirement that $F$ contracts volumes, it is reasonable to hope that H\'enon-like mappings are conjugate to H\'enon mappings.

\bigskip\centerline{19}\bigskip\hrule\bigskip

\centerline{\bf 2. Complex Horseshoes}\bigskip

In this chapter and the next, complex analogs of Smale horseshoes (\cf\ \cite{S}, \cite{Mo}) are defined and analyzed. Using a criterion analogous to the one given by Moser \cite{Mo} in the real case, we
will show that many H\'enon mappings are complex horseshoes. In particular, actual H\'enon mappings (of degree 2) are complex horseshoes when $|c|$ is sufficiently large.

In Figure 1.1, only (a) and (c) appear to be horseshoes. Basically, we would like to say that a horizontal H\'enon-like mapping, $F$, of degree $d$ is a complex horseshoe of degree $d$ if the projections
 $$
 \pi_1: \bigcap_{0\le m\le n}F^{\circ m}(B) \to \C \qquad \text{and} \qquad \pi_2: \bigcap_{0\le m\le n}F^{\circ -m}(B) \to \C 
 $$ 
 are trivial fibrations with fibers disjoint unions of $d^n$ discs.

However, this is not general enough for our purposes. When we investigate the existence of horseshoes in the setting of homoclinic points in chapter 4, the natural domains for these 

\bigskip\centerline{20}\bigskip\hrule\bigskip

\noindent 
mappings become quite ``wiggly'' and the meaning of $\pi_1$ and $\pi_2$ as horizontal and vertical projections outside of the domain of definition, $B$, becomes unclear and the notion of any such projections
becomes unusable. Figure 4.1 in Chapter 4 should convince the reader of this.

We now give a definition with weaker conditions, inspired by Proposition 1.2, which encompasses the H\'enon-like mappings of the previous chapter. Now, instead of requiring $B$ to be an actual
bidisc as in Chapter 1, $B$ may be an embedded bidisc. More precisely, letting $D \subset \C$ be the open unit disc, assume that there is an embedding, $\varphi: \overline D^2 \to \C^2$, which is analytic
on $D^2$ and such that $B = \varphi(D^2)$ and, naturally, $\overline B = \varphi(\overline D^2)$. Set the notation
 $$
 \partial B_H = \varphi(\overline D \times \partial \overline D)\qquad \text{and}\qquad \partial B_V = \varphi(\partial \overline D \times \overline D)  
 $$
 for the horizontal and vertical boundaries of $B$. Also, define horizontal and vertical slices: $H_y = \varphi(D \times \{y\})$ and $V_x = \varphi(\{x\} \times D)$ for all $x, y \in D$.

Consider mappings $F: \overline B \to \C^2$ which are injective and continuous on $\overline B$ and analytic on $B$, and such that either
 $$
 F(\overline B) \cap \partial B_H = \emptyset\qquad \text{and}\qquad \overline B \cap F(\partial B_V) = \emptyset\tag2.1
 $$

\bigskip\centerline{21}\bigskip\hrule\bigskip

\noindent 
or
 $$
 \overline B \cap F(\partial B_H) = \emptyset\qquad \text{and}\qquad F(\overline B) \cap \partial B_V = \emptyset.\tag2.2
 $$

Under these conditions, we have the following:

\proclaim {Lemma 2.1} 
 For all $y \in D$, $\pi_1 \circ \varphi^{-1}: F(H_y) \cap B \to D$ is proper. 
 \endproclaim

\demo {Proof} 
 Consider the preimage, $S$, of a compact subset, $C$, of $D$. Suppose a sequence in $S$ converges to a point $\bold x$ which is not in $S$. The second part of either condition (2.1) or (2.2) above implies
that $\bold x$ is not in $B$ while the first part implies that $\bold x$ is not in $\partial B_H$. Compactness of $C$ implies $\bold x$ is not in $\partial B_V$. 
 \newline\null\hfilll{\bf QED}
 \enddemo

Such a proper mapping has a degree and since the degree is integer-valued and continuous in $y$, this defines a constant, the degree of such a mapping $F$. Now a class of mappings generalizing the
H\'enon-like mappings of Chapter 1 can be defined.

\definition {Definition 2.2} 
 $F: \overline B \to \C^2$ is a {\it quasi-H\'enon-like mapping of degree $d$\/} if there exists a mapping $G: \overline B \to \C^2$ such that

\bigskip\centerline{22}\bigskip\hrule\bigskip

 \roster 
 \item"(1)" Both $F$ and $G$ are injective and continuous on $\overline B$ and analytic on $B$.  
 \item"(2)" $F \circ G = \roman{Id}$ and $G \circ F = \roman{Id}$ where each makes sense. \newline
 Hence, rename $G$ as $F^{-1}$.  
 \item"(3)" Either 
 \item"" {\rm (a)}\quad  $F(\overline B) \cap \partial B_H = \emptyset$ and $\overline B \cap F(\partial B_V) = \emptyset$ or 
 \item"" {\rm (b)}\quad  $\overline B \cap F(\partial B_H) = \emptyset$ and $F(\overline B) \cap \partial B_V = \emptyset$. 
 \endroster
 The degree from Lemma 2.1 is $d \ge 2$.  Moreover, call $F$ either {\it horizontal\/} or {\it vertical\/} according to whether it satisfies (a) or (b), respectively.  
 \enddefinition

\remark {Remark} 
 Note that, as in Definition 1.1 for H\'enon-like mappings, the conditions (a) and (b) of Definition 2.2 are dual in the sense that (a) and (b) are equivalent to 
 \roster 
 \item"" {\rm (a')}\quad  $\overline B \cap F^{-1}(\partial B_H) = \emptyset$ and $F^{-1}(\overline B) \cap \partial B_V = \emptyset$ and 
 \item"" {\rm (b')}\quad  $F^{-1}(\overline B) \cap \partial B_H = \emptyset$ and $\overline B \cap F^{-1}(\partial B_V) = \emptyset$, respectively.
 \endroster 
 \endremark

Of course, Proposition 1.2 gives the following:

\proclaim {Proposition 2.3} 
 H\'enon-like mappings of degree $d$ are quasi-H\'enon-like of degree $d$.
 \endproclaim

Now, using the notion of quasi-H\'enon-like mappings, complex 

\bigskip\centerline{23}\bigskip\hrule\bigskip

\noindent 
horseshoes may be defined as suggested above.

\definition {Definition 2.4} 
 A {\it complex horseshoe of degree $d$\/} is a quasi-H\'enon-like mapping of degree $d$, $F: \overline B \to \C^2$, such that, for all integers $n > 0$, depending on if $F$ is horizontal or vertical, then
either the projections 
 $$
 \pi_1 \circ \varphi^{-1}: \bigcap_{0\le m\le n} F^{\circ m}(B) \to \C\qquad \text{and}\qquad \pi_2 \circ \varphi^{-1}: \bigcap_{0\le m\le n} F^{\circ -m}(B) \to \C\hphantom{,}
 $$ 
 or  
 $$
 \pi_2 \circ \varphi^{-1}: \bigcap_{0\le m\le n} F^{\circ m}(B) \to \C\qquad \text{and}\qquad \pi_1 \circ \varphi^{-1}: \bigcap_{0\le m\le n} F^{\circ -m}(B) \to \C,
 $$ 
 respectively, are trivial fibrations with fibers disjoint unions of $d^n$ discs. 
 \enddefinition

\remark {Remarks} 
 (1) This definition should be compared to real Smale horseshoes, where the corresponding projections are bundles of $d^n$ intervals. In the real definition, hyperbolicity conditions must be imposed. Being
in the complex domain, we get them for free.

(2) As for H\'enon-like mappings, we can define sets $K$, $J$, etc. The set $B$ would ideally be a small neighborhood of $K$. 
 \endremark

\bigskip\centerline{24}\bigskip\hrule\bigskip

In the context of H\'enon-like mappings, the following results will show the close relation between complex horseshoe mappings of degree $d$ and polynomial-like mappings of degree $d$ whose critical points
escape immediately.

The following definition, borrowed from Moser \cite{Mo}, is the key tool in this study of complex horseshoes. Apparently, this concept is due to V. Alekseev \cite{A}.

Let $M$ be a differentiable manifold, $U \subset M$ an open subset, and $f: U \to M$ a differentiable mapping.

\definition {Definition 2.5} 
 A field of cones $\pmb{\Cal C} = (\Cal C_{\bold x} \subset T_{\bold x}M)_{x\in U}$ on $U$ is an {\it $f$-trapping\/} field if
 \roster 
 \item"(1)"  $\Cal C_{\bold x}$ depends continuously on $\bold x$ and
 \item"(2)"  whenever $\bold x \in U$ and $f(\bold x) \in U$, then $d_{\bold x}f(\Cal C_{\bold x}) \subset \Cal C_{f(\bold x)}$. 
 \endroster 
 \enddefinition

Now consider the connection between trapping fields of cones and complex horseshoes.

\proclaim {Theorem 2.6} 
 Let $F: \overline B \to \C^2$ be a quasi-H\'enon-like mapping of degree $d$. The following are equivalent: 
 \roster 
 \item"(1)"  $F: \overline B \to \C^2$ is a complex horseshoe of degree $d$

\bigskip\centerline{25}\bigskip\hrule\bigskip

\noindent 
 \item"(2)"  There exist continuous, positive functions $\alpha({\bold z})$ and $\beta({\bold z})$ on $B$ such that the field of cones
 $$
 \Cal C_{\bold z} = \bigl\{\,(\xi_1, \xi_2)\bigm| |\xi_2| < \alpha ({\bold z})|\xi_1|\,\bigr\}
 $$ 
 is $F$-trapping and the field of cones 
 $$
 \Cal C_{\bold z}' = \bigl\{\,(\xi_1, \xi_2)\bigm| |\xi_1| < \beta ({\bold z})|\xi_2|\,\bigr\}
 $$ 
 is $F^{-1}$-trapping. 
 \item"(3)"  $F(\overline B) \cap \overline B$ and $F^{-1}(\overline B) \cap \overline B$ both have $d$ connected components. 
 \endroster 
 \endproclaim

\remark {Remarks} 
 Note that (2) $\Rightarrow$ (1) is borrowed from Moser and that the implication (1) $\Rightarrow$ (2) is what the contractive nature of complex analytic mappings gives us for free. (3) arises naturally in
the proof of (2) $\Rightarrow$ (1). When considering a mapping which is actually H\'enon-like, consideration of the critical points of $F_{1,y}$ and $F_{2,x}^{-1}$ (or $F_{2,x}$ and $F_{1,y}^{-1}$) in
light of the equivalences above yields the following. 
 \endremark

\proclaim {Corollary 2.7} 
 Let $F: \overline{D_1}\times \overline{D_2} \to \C^2$ be a H\'enon-like mapping of degree $d$. The following are equivalent: 
 \roster 
 \item"(1)"  $F: \overline{D_1} \times \overline{D_2} \to \C^2$ is a complex horseshoe of degree $d$. 
 \item"(2)"  For all $(x, y) \in \overline{D_1} \times \overline{D_2}$, the critical values of the 

\bigskip\centerline{26}\bigskip\hrule\bigskip

\noindent 
polynomial-like mappings $F_{1,y}$ and $F_{2,x}^{-1}$ (or $F_{2,x}$ and $F_{1,y}^{-1}$) lie outside of
$\overline{D_1}$ and $\overline{D_2}$, respectively. 
 \endroster
 \endproclaim

\demo {Proof of Theorem 2.6} 
 Without loss of generality we can assume that $B = D^2$ and that $F$ satisfies condition (a) of Definition 2.2 [or Definition 1.1 in case it is actually H\'enon-like].

First we shall prove that (1) $\Rightarrow$ (2). Let $F: \overline D^2 \to \C^2$ be a complex horseshoe of degree $d$, let ${\bold z} = (x, y) \in D^2$ and set
 $$
 \Cal C_{\bold z} = \bigl\{\,(\xi_1,\xi_2)\bigm| |\xi_2 |_y < |\xi_1|_x\,\bigr\},
 $$ 
 where $|\xi|_x$ means the Poincar\'e length at $x$, \ie, 
 $$
 |\xi|_x = |\xi|/\left(1 - |x|^2\right),
 $$ 
 so that this field of cones is defined by 
 $$
 \alpha({\bold z}) = \frac{1 - |y|^2}{1 - |x|^2}.
 $$ 

We will show that this field is $F$-trapping. First note that (1) $\Rightarrow$ (3). $F^{-1}(D^2) \cap D^2$ has $d$ connected components---call them $U_1$, $U_2$, \dots, $U_d$---and $F(D^2) \cap D^2$
also has $d$ connected 

\bigskip\centerline{27}\bigskip\hrule\bigskip

\noindent 
components, $F(U_1)$, $F(U_2)$, \dots, $F(U_d)$. Moreover, $F: U_i \to F(U_i)$ is an analytic isomorphism for each $i = 1, \dots, d$.

Now suppose that ${\bold x} \in D^2$, $F({\bold x}) \in D^2$, and $\xi \in \Cal C_{\bold x} \subset T_{\bold x}D^2$. Take an analytic $\alpha: D \to D$ with $\xi$ tangent to the graph of $\alpha$,
which we write as $\text{gr} (\alpha)$, at ${\bold x}$. Suppose that ${\bold x} \in \text{gr} (\alpha) \cap U_i$ and consider $F(\text{gr} (a) \cap U_i)$. Define $\beta: D \to D$ such that $\text{gr}
(\beta) = F(\text{gr} (\alpha) \cap U_i)$. So $d_{\bold x}F(\xi)$ is tangent to the graph of $\beta$. Since $F(\text{gr}(\alpha))$ has $d$ components, $\beta$ is not surjective. By Schwarz' Lemma, $\beta$ is
contracting in the Poincar\'e metric and $d_{\bold x}F(\xi) \in \Cal C_F({\bold x})$. Thus, (1) $\Rightarrow$ (2).

We next show that (2) $\Rightarrow$ (3). Suppose that $\Cal C_{\bold z}$ and $\Cal C_{\bold z}'$ are $F$-trapping and $F^{-1}$-trapping fields of cones, respectively. Consider $\pi_1
\circ F: H_y \to \C$ for an arbitrary $y \in D$. The critical values of this mapping are outside $D$, since otherwise the critical points ${\bold z_i} \in H_y$ and their images would lie in $D^2$, and the
image of horizontal vectors at ${\bold z_i}$ would be vertical vectors at $F({\bold z_i})$, and this is forbidden by the existence of the trapping field. By Lemma 2.2 and the definition of
quasi-H\'enon-like of degree $d$, $F(\overline B) \cap \overline B$ has $d$ components. This proves that (2) $\Rightarrow$ (3).

\bigskip\centerline{28}\bigskip\hrule\bigskip

\noindent 
It remains to show that (3) $\Rightarrow$ (1). So, assume that $F(\overline B) \cap \overline B$ has $d$ components. Suppose by induction that, for all $n = 1, 2, \dots, N - 1$,
 $$
 \pi_2: \bigcap_{0\le m\le n} F^{\circ -m}(D^2) \to D \tag{1}
 $$ 
 is a trivial fibration with fiber $d^n$ discs and \newline
 (2) for all $y \in D$ and each component $U$ of 
 $$
 H_y \cap \bigcap_{0\le m\le n} F^{\circ -m}(D^2)
 $$ 
 the map $\pi_1 \circ F^{\circ n}: U \to D$ is an analytic isomorphism.

Let us prove it for $n = N$. Choose such a component $U$ for $n = N - 1$, and consider $F^{\circ N}(U)$. Since $\pi_1 \circ F^{\circ N-1}$ is an isomorphism from $U$ to $D$ and $F(\overline D^2) \cap
\overline D^2$ has $d$ components, $F^{\circ N}(U) \cap \overline D^2$ has $d$ components also, $U_1$, $U_2$, \dots, and $U_d$, each homeomorphic to a disc, and for each the mapping $\pi_1: U_i \to D$ is an
analytic isomorphism. So, (3) $\Rightarrow$ (1). 
 \newline\null\hfilll{\bf QED}
 \enddemo

\bigskip\centerline{29}\bigskip\hrule\bigskip
 
\demo {Proof of Corollary 2.7} 
 Since polynomial-like of degree $d$ with critical points escaping implies an image with $d$ components, (2) implies that $F(\overline D^2) \cap \overline D^2$ has $d$ components. So, (2) $\Rightarrow$
(1) follows from (3) $\Rightarrow$ (1) of Theorem 2.6.

By (1) $\Rightarrow$ (2) of Theorem 2.6, there exist $F$-trapping and $F^{-1}$-trapping fields of cones. Consider $F_{1,y} = \pi_1 \circ F: H_y \to \C$ for an arbitrary $y \in D$. This is polynomial-like of
degree $d$. The critical values of this mapping are outside $D$, since otherwise the critical points ${\bold z_i} \in H_y$ and their images would lie in $D^2$, and the image of horizontal vectors at ${\bold
z_i}$ would be vertical vectors at $F({\bold z_i})$, and this is forbidden by the existence of the trapping field. Therefore, (1) $\Rightarrow$ (2). 
 \newline\null\hfilll{\bf QED}
 \enddemo

\remark {Remarks} 
 The definition of a trapping field of cones made above is designed to be simple and easily verifiable, but is not the strongest definition which will still give the result above. Such a definition would
require the field of cones to be defined only on the set $\cap_{n\ge0} f^{\circ n}(U)$ and could be used to obtain sharper results (\cf Proposition 2.10 and the remarks following it).
 \endremark

\proclaim {Proposition 2.8} 
 The diameters of the discs in the fibers 

\bigskip\centerline{30}\bigskip\hrule\bigskip

\noindent 
above tend to $0$ with $n$. 
 \endproclaim

The essential part of the proof of this proposition is the following lemma on pol\-y\-nom\-i\-al-like mappings. Let $D$ be an open disc and let $\alpha: D \to D$ be an analytic injection. The mapping $f =
f_\alpha: D \to \C$ defined by $f(x) = f_\alpha(x) = \pi_1(F(x, \alpha_j(x))$ is polynomial-like of degree $d$ with the critical points escaping immediately, \ie, with the critical values in $\CmD$. So,
there exist $d$ components, $U_1$, $U_2$, \dots, and $U_d$, of $f^{-1}(D)$ and $d$ mappings, $g_i = g_{\alpha,i}: D \to U_i$, which are inverses of $f|_{U_i}$, for $i = 1, 2, \dots, d$.

\proclaim {Lemma 2.9} 
 Independent of $\alpha$ and $i$, the diameter of the image of a connected subset $U$ of $D$ under $g_i$ as compared with the diameter of $U$ shrinks by a constant factor less than one. 
 \endproclaim

\demo {Proof} 
 By the definition of polynomial-like mappings it follows that the distance between the boundary of $D$ and the boundary of $g_i(D)$ is bounded away from $0$. Now, the distance between two points, $x$ and
$y$, of $g_i(D)$ equals the distance between the points $f(x)$ and $f(y)$ when measured in the Poincar\'e metrics on $g_i(D)$ and $D$, respectively. Since $g_i(D) \subset D$, when 

\bigskip\centerline{31}\bigskip\hrule\bigskip

\noindent 
measured in $D$, the distance between $x$ and $y$ is always less than the distance between $f(x)$ and $f(y)$. There exists a constant $K < 1$ such that $d_D(x, y) \le K{\,}d_D(f(x), f(y))$. So, $K$ is the
shrinking factor. 
 \newline\null\hfilll{\bf QED}
 \enddemo

\demo {Proof of Proposition 2.8} 
 Consider a component, $V$, of a fiber over $y$ of  $\pi_2: \underset{0\le m\le n}\to\cap F^{\circ -m}(D^2) \to D$. We have 
 $$
 V = V_n \subset V_{n-1} \subset V_{n-2} \subset \dotsb
 $$ 
  where each $V_j$ is the component of $H_y \cap F^{\circ -j}(D^2)$ which contains $V$. Each $V_j$ has the property that $F^{\circ j}(V_j)$ is the graph of an analytic injection $\alpha_j: D \to D$. Now
there exist such mappings $\alpha_1$, $\alpha_2$, \dots, $\alpha_n: D \to D$ and integers $i_1$, $i_2$, \dots, $i_n$, which are between $1$ and $d$, inclusive, such that
 $$
 \pi_2(V) = g_{\alpha_n,i_n} \circ \dots \circ g_{\alpha_2,i_2} \circ g_{\alpha_1,i_1}(D).
$$ 
By Lemma 2.9, each $g_{\alpha_j,i_j}$ contracts the diameter of $D$ by a
constant factor $K < 1$ independent of $\alpha_j$ and $i_j$. 
 \newline\null\hfilll{\bf QED}
 \enddemo

This criterion can be used to show that for each $a$ there 

\bigskip\centerline{32}\bigskip\hrule\bigskip

\noindent 
exists $r(a)$ such that if $|c| > r(a)$, then the H\'enon mapping $F_{a,c}$ is a complex horseshoe. Of course, in the real locus this was known
(\cf\ \cite{DN}, \cite{N}), except that then $c$ has to be taken very negative. When $c$ is large and positive, all the ``horseshoe behavior'' is complex. 

\proclaim {Proposition 2.10} 
 For each $a \ne 0$ and each $c$ such that $|c| > \left(5/4 + \sqrt5/2\right)(1 + |a|)^2$,  there exists an $R$ such that $F_{a,c}: \overline{D_R}^2 \to \C^2$ is a complex horseshoe.
 \endproclaim

\demo {Proof} 
 Set $F = F_{a,c}$ and set 
 $$
 R = \alpha\left(1 + |a| + \sqrt{(1 + |a|)^2 + 4|c|}\right)\tag2.3
 $$ 
 with $\alpha > 1/2$. By Proposition 1.4 and example (1) which follows it, $F: \overline{D_R}^2 \to \C^2$ is quasi-H\'enon-like. By Theorem 2.6, this proposition can be proved by
showing that, for appropriate values of $a$, $c$, and $\alpha$, the constant field of cones defined simply by 
 $$
 \bigl\{\,(\xi_1, \xi_2)\bigm| |\xi_2| < |\xi_1|\,\bigr\}
 $$ 
 is $F$-trapping, and, similarly, that the field of cones
 $$
 \bigl\{\,(\xi_1, \xi_2)\bigm| |\xi_1| < |\xi_2|\,\bigr\}
 $$ 

\bigskip\centerline{33}\bigskip\hrule\bigskip

\noindent 
 is $F^{-1}$-trapping.

The key to showing that the former field of cones is $F$-trapping is the observation that $F(x,y) \in \overline{D_R}^2$ implies that $|x^2 + c - ay| \le R$ which implies that 
 $$
 |x|^2 \ge |c| - R(1 + |a|).\tag2.4
 $$ 
 Since $d_{(x,y)}F(\xi_1, \xi_2) = (2x\xi_1 - a\xi_2, \xi_1)$, we want to show that
 $$
 |\xi_1| > |\xi_2|\qquad \text{implies}\quad |2x\xi_1 - a\xi_2| > |\xi_1|.
 $$
 Using \thetag{2.4} to substitute for $x$, it suffices to show that 
 $$
 2\sqrt{|c| - R(1 + |a|)} - |a| > 1.\tag2.5
 $$ 
 A similar analysis for $F^{-1}$ and the latter field of cones yields the same inequality \thetag{2.5}.

After using (2.3) to substitute into (2.5) for $R$ and setting $|c| > \beta(1 + |a|)^2$, a computation gives
 $$
 |c| > (\alpha + \tfrac14 + \alpha\sqrt{1 + 4\beta}) (1 + |a|)^2,
 $$

\bigskip\centerline{34}\bigskip\hrule\bigskip

\noindent 
which is satisfied when 
 $$
 \beta \ge 2\alpha^2 + \alpha + \tfrac14 + \alpha\sqrt{4\alpha^2 + 4\alpha + 2}.
 $$ 
 As $\alpha$ goes to $\frac12$, the lower bound for $\beta$ goes to $5/4 + \sqrt5/2$. 
 \newline\null\hfilll{\bf QED}
 \enddemo

\remark {Remarks}
(1) The above result says essentially everything about H\'enon mappings in the parameter range to
which it applies.

(2) In the next section, it is shown that a complex horseshoe of degree $d$, restricted to the set $K$, is conjugate to the full shift on $d$ symbols. Milnor \cite{M} has used a completely different method
to show the existence of an embedding of the shift for $\beta = 2$. The degenerate case, $a = 0$, shows that this is the strongest result given by the condition $|c| > \beta(1 + |a|)^2$ which we could  hope
for.

(3) As was indicated in the remarks after the proof of Theorem 2.6, if we strengthen the definition of trapping fields, then we can improve Proposition 2.10 so that we have complex horseshoes whenever $|c|
> \beta(1 + |a|)^2$ for some $\beta$ between $2$ and $5/4 + \sqrt5/2 \approx 2.368$. By requiring that the trapping field be defined only on the set of points which remain in the bidisc for 

\bigskip\centerline{35}\bigskip\hrule\bigskip

\noindent 
additional iterations, we could get a sequence of improvements for the lower bound of $\beta$. However, it was intended that the existence of trapping fields be easily verified and hence they were defined
on the bidisc rather than on $K_+$ or some compromise between these two. Nevertheless, with much inequality manipulating sharper results can be achieved.

(4) Allowing fields of larger cones such as 
$$
\bigl\{\,(\xi_1, \xi_2)\bigm| \gamma|\xi_2| < |\xi_1|\,\bigr\},
$$ 
with $0 < \gamma < 1$, it is possible to get better estimates for lower bounds depending on $|a|$ for $|c|$ so that $F$ is a complex horseshoe. For large $|a|$ it is not difficult to get sharper estimates
than those given above. For example, for $|a| = 1$, it is possible to use very large cones to see that horseshoes exist for $|c| > 7.2$, as compared with $|c| > 9.48$ using the result above (or $|c| > 8$
using Milnor's result). 

(5) Using Example (2) following Proposition 1.4 and the same procedure as in the proof of Proposition 2.10 above, we can get examples of complex horseshoes of higher degrees. Of course, sharper results are
possible here also.
 \endremark

\bigskip\centerline{36}\bigskip\hrule\bigskip

\centerline{\bf 3. Complex Horseshoes as Shift Dynamical Systems}\bigskip

Suppose $d$ is an integer greater than or equal to two. Let $\Cal S_d = \{\,0, 1, \dots, d\,\}^{\Z}$ and $\tau_d: \Cal S_d \to \Cal S_d$ be the shift on $d$ symbols. Let $\Z_+$ be the positive integers and
$\Z_- = \Z \setminus\Z_+$. Let 
 $$
 \Cal S_d^+ = \{\,0, 1, \dots, d\,\}^{\Z_+} \qquad\text{and} \qquad\Cal S_d^- = \{\,0, 1, \dots, d\,\}^{\Z_-}.
 $$
 
Let $\tau_d^+: \Cal S_d^+ \to \Cal S_d^+$ and $\tau_d^-: \Cal S_d^- \to \Cal S_d^-$ be the corresponding one-sided shifts.

Let $F: \overline D^2 \to \C^2$ be a complex horseshoe of degree $d$. The following results allow a complete understanding of the invariant subsets for a complex horseshoe, essentially identical with what
you get for real horseshoes.

\proclaim {Theorem 3.1} 
 There exists a homeomorphism $\Phi: K \to \Cal S_d$ which conjugates $F$ to the full shift $\tau_d$. 
 \endproclaim

\remark {Remark} 
 In the case where $d = 2$, it will be seen that $\Phi$ is unique up to the automorphism of $\Cal S_2$ exchanging $0$ and $1$. 
 \endremark

In the following we will fix $d$ and set $\Cal S = \Cal S_d$, $\Cal S_+ = \Cal S_d^+$, 

\bigskip\centerline{37}\bigskip\hrule\bigskip

\noindent 
$\Cal S_- = \Cal S_d^-$, $\tau = \tau_d$, $\tau_+ = \tau_d^+$, and $\tau_- = \tau_d^-$. In fact, assume that $d = 2$ without loss of generality.

Theorem 3.1 is actually a corollary of the following result:

\proclaim {Theorem 3.2}  There exist homeomorphisms 
 $$
 \Gamma_+: K_+ \to \Cal S_+ \times D\qquad and\qquad \Gamma_-: K_- \to \Cal S_- \times D
 $$
 which can be written as 
 $$
 \Gamma_+(x, y) = (\gamma_+(x, y), y)\qquad and\qquad \Gamma_-(x, y) = (x, \gamma_-(x, y)),
 $$
 such that the following diagrams commute:
 $$
 \CD 
 K_+ \cap F^{-1}(D^2) @>>\gamma_+> \Cal S_+ @. \qquad @. \qquad @. \qquad @. K_- \cap F(D^2) @>>\gamma_-> \Cal S_- \\
 @VVFV @VV\tau_+V @.\quad\text{and} @. @. @VVF^{-1}V @VV\tau_-V \\ 
 K_+ @>>\gamma_+> \Cal S_+ @. \qquad @. \qquad @. \qquad @. K_- @>>\gamma_-> \Cal S_- 
 \endCD
 $$
 The inclusions $K \subset K_+$ and $K \subset K_-$ induce the canonical projections $\Cal S \to \Cal S_+$ and $\Cal S \to \Cal S_-$. 
 \endproclaim

We will define in Lemma 3.3 models for $\Cal S_+$ and $\Cal S_-$ which are defined which are more closely related to the dynamics of the horseshoe and then

\bigskip\centerline{38}\bigskip\hrule\bigskip

\noindent 
restate the result as Theorem 3.4 in a form which uses this model and such that the new statement is strictly stronger than Theorem 3.2.

Choose a polynomial-like mapping $f: U' \to U$ of degree $2$ such that the critical value lies in $U \setminus \overline{U'}$ (for instance, $U' = D_3(0)$, $U = D_9(4)$, and $f(z) = z^2 + 4$). Since the
critical value is in $U \setminus \overline{U'}$, the set $f^{-1}(U')$ consists of two open sets, $U_0$ and $U_1$, homeomorphic to discs, and each $f|_{U_i}: U_i \to U'$ is an isomorphism. Let $g_i: U' \to
U_i$ be the corresponding inverses.

\proclaim {Lemma 3.3} 
 \roster  
 \runinitem"(a)" For each sequence $\varepsilon = (\varepsilon_1, \varepsilon_2, \dots)$ of $0$' s and $1$' s, the nested intersection
 $$
 p_+(\varepsilon) = \bigcap_{n \ge 0} g_{\varepsilon_n} g_{\varepsilon_{n-1}}\dots g_{\varepsilon_1}(U')
 $$
 \item"" reduces to a single point.
 \item"(b)" The mapping $p_+: \Cal S_+ \to U$ is a homeomorphism onto $J_f$, the Julia set of $f$, which conjugates $f$ to the one-sided shift $t_+: \Cal S_+ \to \Cal S_+$. 
 \item"(c)" There is a similar homeomorphism $p_-: \Cal S_- \to J_f$ 

\bigskip\centerline{39}\bigskip\hrule\bigskip

\noindent 
which conjugates $f$ to the one-sided shift $\tau_-: \Cal S_- \to \Cal S_-$.
 \endroster 
 \endproclaim

\demo {Proof} 
 Part (a) is proved by contraction as in Proposition 2.8 and is simpler. Actually, the proof of Proposition 2.8 was adapted from this lemma. Part (b) follows from the fact that $J_f = K_f$ when the critical
value is not in $U'$.
 \newline\null\hfilll{\bf QED}
 \enddemo

\proclaim {Theorem 3.4} 
 There exist homeomorphisms
 $$
 \align 
 \Psi_+: D^2 \to U \times D\qquad &\text{and}\qquad \Psi_-: D^2 \to U \times D \\ 
 \intertext{which can be written as} 
 \Psi_+(x, y) = (\psi_+(x, y), y)\qquad &\text{and}\qquad \Psi_-(x, y) = (x, \psi_-(x, y))
 \endalign
 $$
 such that the following diagrams commute:
 $$
 \CD 
 D^2 \cap F^{-1}(D^2) @>>\psi_+> U' @. \qquad @. \qquad @. \qquad @. D^2 \cap F(D^2) @>>\psi_-> U' \\ 
 @VVFV @VVfV @.\quad\text{and} @. @. @VVF^{-1}V @VVfV \\ 
 D^2 @>>\psi_+> U @. \qquad @. \qquad @. \qquad @. D^2 @>>\psi_-> U 
 \endCD
 $$
 \endproclaim

\demo {Proof} 
 First choose the restriction

\bigskip\centerline{40}\bigskip\hrule\bigskip

 $$
 \psi_+: \partial D \times D \to \partial U
 $$
 to be induced from a homeomorphism $\partial D \to \partial U$. Next label the two components $V_0$ and $V_1$ of $D^2 \cap F^{-1}(D^2)$, and define $\psi_+$ on $\partial V_i$ to be $g_i \circ
\psi_+ \circ F$. 
 \enddemo

\proclaim {Lemma 3.5} 
 The restricted mapping
 $$
 \Psi_+: \partial V_i \to \partial U_i \times D \qquad \text{defined by} \qquad \Psi_+(x, y) = (\psi_+(x, y), y)
 $$
 is a homeomorphism. 
 \endproclaim

\demo {Proof} 
 It is enough to prove that for each $y \in D$, the restriction
 $$
 \psi_+: \partial V_i \cap D_y \to \partial U_i
 $$
 is a homeomorphism. It is of degree $1$, and locally injective because $F(\partial V_i \cap D_y)$ is a simple closed curve with tangent vectors in the trapping field, hence transversal to the (vertical)
fibers of $\psi_+$ on $\partial D \times D$. But an immersion of degree $1$ from a simple closed curve to a simple closed curve is a
homeomorphism. 
 \newline\null\hfilll{\bf QED for Lemma 3.5}
 \enddemo

\bigskip\centerline{41}\bigskip\hrule\bigskip

The following is a simple result on surface topology.

\proclaim {Lemma 3.6} 
 If $X_1$ and $X_2$ are spheres with three holes and $h: \partial X_1 \to \partial X_2$ is an orientation-preserving homeomorphism, then there exists a homeomorphism $\tilde h: X_1 \to X_2$ extending $h$. 
 \endproclaim

\demo {Completion of Proof of Theorem 3.4} 
 We can extend $\psi_+$ to give a homeomorphism
 $$
 \Psi_+: D^2 \setminus (V_0 \cup V_1) \to ( U \setminus (U_0 \cup U_1) ) \times D
 $$
 by the Cerf fibration theorem.

Let $W = D^2 \setminus (V_0 \cup V_1)$. Now $W \cap F^{-1}(W)$ has two components, $W_0 \subset V_0$ and $W_1 \subset V_1$. Extend $\psi_+$ to $W_i$ by $\psi_+(x, y) = g_i\bigl(\psi_+(F(x, y))\bigr)$. Again
we must prove that $(x, y) \mapsto (\psi_+(x, y), y)$ is a homeomorphism $W_i \to U \setminus (U_0 \cup U_1) \times D$; the proof is just as before.

Continue to extend in exactly the same way to $W \cap F^{\circ -2}(W)$, $W \cap F^{\circ -3}(W)$, \dots, verifying at each stage that the extension is a homeomorphism.

\bigskip\centerline{42}\bigskip\hrule\bigskip

When we are done, we will have a mapping
 $$
 \psi_+: D^2 \setminus K_+ \to U \setminus J_f
 $$
 such that the mapping
 $$
 \Psi_+: D^2 \setminus K_+ \to (U \setminus J_f) \times D\qquad \text{defined by}\qquad \Psi_+(x, y) = (\psi_+(x, y), y)
 $$
 is a homeomorphism and
 $$
 \psi_+ \circ F(x, y) = f \circ \psi_+(x, y)\qquad \text{on}\qquad (D^2 \cap F^{-1}(D^2)) \setminus K_+.
 $$
 In each slice $D_y$, the set $D_y \cap K_+$ is a Cantor set, so that $\overline{D_y} \cap K_+$ is the set of ends of $\overline{D_y} \setminus K_+$. In the same way, $J_f$ is the set of ends of $\overline
U \setminus J_f$. Since the end-point compactification is functorial, the map $\psi_+ extends to D^2$, and since the extension is a homeomorphism $D_y \to U$ for each $y$, the map $\Psi_+$ is a
homeomorphism. 
 \newline\null\hfilll{\bf QED for Theorem 3.4}
 \enddemo

\remark {Remark} 
 It looks as if there were two independent choices of labeling $V_0$ and $V_1$, the two components of $D^2 \cap F^{-1}(D^2)$, and $V_0'$ and $V_1'$, the two components of $D^2 \cap F(D^2)$. However,
 $$
 F: D^2 \cap F^{-1}(D^2) \to D^2 \cap F(D^2)
 $$

\bigskip\centerline{43}\bigskip\hrule\bigskip

\noindent is a homeomorphism. So choose the labelings so that $F(V_0) = V_0'$.

Now Theorem 3.2 follows immediately from the strictly stronger Theorem 3.4. We return to the theorem on conjugating $F$ on $K$ to the full shift on $2$ symbols and use Theorem 3.2.
 \endremark

\demo {Proof of Theorem 3.1} 
 Now the mapping given by $(x, y) \mapsto (\gamma_-(x, y) , \gamma_+(x, y))$ clearly has all the properties required. 
 \newline\null\hfilll{\bf QED for Theorem 3.1}
 \enddemo

We will now show that the dynamics of the projective limit under $f$ underlies the dynamics of the horseshoe mapping. Define $U'_- = \varprojlim (U', f)$ and let $\pi: U_- \to U$ be the projection
defined by $(\dots, z_2, z_1, z_0) \mapsto z_0$. Also, let $U'_- = \pi^1(U') \subset U_-$.

\proclaim {Theorem 3.7}
 There exist homeomorphisms $\eta_\pm: K_\pm \to U_-$ such that the following diagrams commute:
 $$
 \CD 
 K_+ \cap F^{-1}(D^2) @>>\eta_+> U'_- @. \qquad @. \qquad @. \qquad @. K_- \cap F(D^2) @>>\eta_-> U'_- \\ 
 @VVFV @VVP_-V @.\quad\text{and} @. @. @VVFV @VVP_-V \\ 
 K_+ @>>\eta_+> U_- @. \qquad @. \qquad @. \qquad @. K_- @>>\eta_-> U_- 
 \endCD
 $$
 \endproclaim

\bigskip\centerline{44}\bigskip\hrule\bigskip

\demo {Proof} 
 Define $\eta_+$ by
 $$
 \eta_+(x, y) = (\dots, \psi_-(F^{\circ 2}(x, y)), \psi_-(F(x, y)), \psi_-(x, y) ).
 $$
 \newline\null\hfilll{\bf QED}
 \enddemo

\bigskip\centerline{45}\bigskip\hrule\bigskip

\centerline{\bf 4. The Ubiquity of Complex Horseshoes}\bigskip

Suppose that $F: \C^2 \to \C^2$ is an analytic mapping. A point $\bold q$ is a {\it homoclinic point\/} of $F$ if there exists a positive integer $k$ such that
 $$
 \lim_{n \to \infty}F^{\circ kn}(\bold q) = \lim_{n \to \infty}F^{\circ -kn}(\bold q) = \bold p\tag{$*$}
 $$ 
 where, in particular, the limits exist.

Note that the limit point, $\bold p$, in ($*$) is a hyperbolic periodic point of $F$ of period least $k$ satisfying ($*$). To rephrase this definition, $\bold q$ is in both the stable and unstable manifolds
of $F$ at $\bold p$, which will be denoted by $W^s$ and $W^u$, respectively. We will call a homoclinic point {\it transversal\/} if these invariant manifolds intersect transversally.

Note that the invariant manifolds $W^s$ and $W^u$ tend to intersect in lots of points in $\C^2$ unless $F$ is linear or affine. This is quite different from the case of these invariant manifolds in $\R^2$.
These intersections are almost always transversal.

\bigskip\centerline{46}\bigskip\hrule\bigskip

In the real, Smale \cite{S} (\cf\ \cite{Mo}) showed that in a neighborhood of a transversal homoclinic point there exist horseshoes. In this chapter we give an analogous result in the case of complex
horseshoes.

Throughout this chapter we will assume that $F$ is analytic, $\bold p$ is a hyperbolic periodic point of period $k$, and $\bold q$ is a transversal homoclinic point as described above.

\proclaim {Theorem 4.1} 
 For every positive integer $d \ge 2$, there exists an embedded bidisc, $B_d$, centered at $\bold p$ and a positive integer $N = N(d)$ such that $F^{\circ kN}: B_d \to \C^2$ is a complex horseshoe of
degree $d$. 
 \endproclaim

\remark {Remark} 
 By Theorem 3.1, there exist homeomorphisms $\Phi_d: B_d \to \Cal S_d$ which conjugate $F^{\circ kN(d)}$ to the full shifts $\tau_d$ on $d$ symbols.

Without loss of generality and in order to simplify the notation, assume that $k = 1$, \ie, that $\bold p$ is a fixed point. Also, fix a positive integer $d \ge 2$.

Essentially, the proof of Theorem 4.1 involves two steps. 

\bigskip\centerline{47}\bigskip\hrule\bigskip

\noindent First, it must be shown that there are quasi-H\'enon-like mappings floating around and, second, that these yield complex horseshoes. We proceed by doing the essential part of the second step first,
that is, showing the existence of trapping fields of cones. 
 \endremark

\proclaim {Lemma 4.2} 
 There exists a neighborhood $U$ of $\bold p$ such that the field of cones 
 $$
 \Cal C_{\bold z} = \bigl\{\,(\xi_1, \xi_2)\bigm| |\xi_2| < |\xi_1|\,\bigr\}
 $$ 
 is $F$-trapping on $U$ and the field of cones 
 $$
 \Cal C_{\bold z}' = \bigl\{\,(\xi_1, \xi_2)\bigm| |\xi_1| < |\xi_2|\,\bigr\}
 $$ 
 is $F^{-1}$-trapping on $U.$ 
 \endproclaim

\demo {Proof} 
 Suppose that the eigenvalues of $d_{\bold p}F$ are $\mu$ and  $\lambda$ with $0 < |\mu| < 1 < |\lambda|$. In a small neighborhood, $U$, of $\bold p$, we can use a set of local coordinates $(u, v)$ such
that 
 $$
 F(u, v) = \bigl(f(u, v), g(u, v)\bigr)
 $$ 
 with the properties that 
 $$
 f(0, v) = g(u, 0) = 0
 $$

\bigskip\centerline{48}\bigskip\hrule\bigskip

\noindent and the partial derivatives 
 $$
 f_u(0, 0) = \lambda \qquad \text{and} \qquad f_v(0, 0) = \mu.
 $$ 
 Consider an arbitrary point $(u_0, v_0) \in U$ and fix a tangent vector $(\xi_0, \eta_0) \in \Cal C_{(u_0, v_0)}$, that is, such that $|\eta_0 | < |\xi_0|$. Define 
 $$
 (u_1, v_1) = F(u_0, v_0) \qquad \text{and} \qquad (\xi_1, \eta_1)  =  d_{(u_0, v_0)}F(\xi_0, \eta_0).
 $$ 
 Assume that $U$ is a small bidisc $D_r^2$ in $(u, v)$ coordinates. Now compute: 
 $$
 \align 
 |\xi_1| &= |f_u(\xi_0) + f_v(\eta_0)|\\ 
 &\ge |f_u(\xi_0)| - |f_v(\eta_0)|\\ 
 &\ge |\lambda| |\xi_0| - O(r)|\eta_0|\\ 
 &\ge \bigl(|\lambda| - O(r)\bigr)|\xi_0|\\ 
 \intertext{and}
 |\eta_1| &= |g_u(\xi_0) + g_v(\eta_0)|\\ 
 &\le |g_u(\xi_0)| + |g_v(\eta_0)|\\ 
 &\le O(r) |\xi_0| + |\mu | |\eta_0|\\ 
 &\le \bigl(|\mu| + O(r)\bigr)|\xi_0|.
 \endalign
 $$ 

\bigskip\centerline{49}\bigskip\hrule\bigskip

\noindent So, for $r$ small enough, $|\xi_1| > |\eta_1|$ and $d_{(u_0, v_0)}F(\xi_0, \eta_0) \in \Cal C_{F(u_0, v_0)}$. 
 \newline\null\hfilll{\bf QED}
 \enddemo

Now, we show that some iterate of $F$ is quasi-H\'enon-like of degree $d$ on a bidisc contained in $U$.

\bigskip
 \centerline{\epsfbox{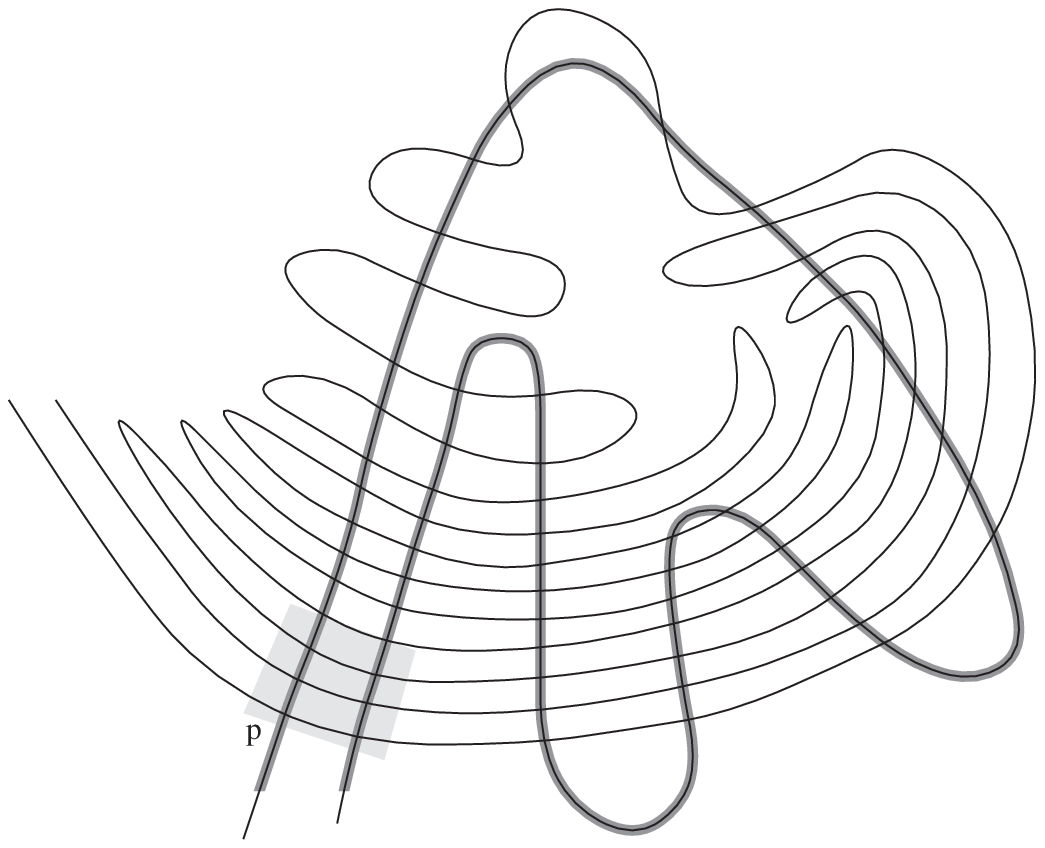}}
 \smallskip
 \centerline{Figure 4.1:{\quad}Horseshoes from transverse homoclinic points}
 \medskip

\bigskip\centerline{50}\bigskip\hrule\bigskip

\proclaim {Lemma 4.3}
 There exist $D_s \subset W^s$ and $D_u \subset W^u$ isomorphic to discs with $D_s, D_u \subset U$ and a positive integer $n$ such that $F^{\circ n}(D_u)$ intersects $D_s$ in exactly $d$ points and the
following conditions hold: 
 $$
 F^{\circ n}(D_u) \cap \partial \overline{D_s} = \emptyset \qquad \text{and} \qquad F^{\circ n}(\partial \overline{D_u}) \cap D_s  = \emptyset 
 $$ 
 \endproclaim

\demo {Proof} 
 Choose subsets $D_s \subset W^s \cap U$ and $D_u \subset W^u \cap U$ such that $D_s \cap D_u = \{\bold p\}$ and both $D_s$ and $D_u$ are isomorphic to discs. Since $W^s = \cup_{j>0} F^{\circ j}(D_u)$
and $W^s$ and $D_u$ intersect in homoclinic points near $\bold p$, there exists a smallest positive integer $n$ such that $F^{\circ n}(D_u) \cap D_s$ contains at least $d$ points, including $\bold p$.

Note that the set $F^{\circ n}(D_u) \cap D_s$ must consist of finitely many points for all positive integers $n$.

If $F^{\circ n}(D_u) \cap \partial \overline{D_s} \ne \emptyset$, then take $D_s$ to be slightly smaller. If $F^{\circ n}(\partial \overline{D_u}) \cap D_s \ne \emptyset$, then take $D_u$ to be slightly
smaller. In either case, make the adjustments to keep the same number of points in $F^{\circ n}(D_u) \cap D_s$.

If there are more than $d$ points in $F^{\circ n}(D_u) \cap D_s$, then deform $D_u$ slightly, making it smaller, to exclude some points from $F^{\circ n}(D_u) \cap D_s$. This can be done in an orderly way.
In 

\bigskip\centerline{51}\bigskip\hrule\bigskip

\noindent particular, there is a natural metric on $D_s$---the Poincar\'e metric of a disc and consider the point $\bold x$ of $F^{\circ n}(D_u) \cap D_s$ which is furthest from $\bold p$ in this metric (or
one such point if more than one has this property). Now deform $D_u$ by taking out $F^{\circ -n}(\bold x)$ and staying clear of the other preimages of points in $F^{\circ n}(D_u) \cap D_s$. 
 \newline\null\hfilll{\bf QED}
 \enddemo

\proclaim {Lemma 4.4} 
 There exists a nonnegative integer $m$ such that, if $D_{u, m} = F^{\circ -m}(D_u)$ and $U_m = D_{u, m} \times D_s$, then $F^{\circ n+m}|_{U_m}$ is a H\'enon-like mapping of degree $d$. 
 \endproclaim

\demo {Proof} 
 For nonnegative integers $j$, set $D_{u, j} = F^{\circ -j}(D_u) \subset D_u$ and $U_m = D_{u, m} \times D_s$. Let $m$ be the smallest nonnegative integer such that $U_m \cap F^{\circ n+m}(\partial D_{u, m}
\times D_s) = \emptyset$. Taking $m$ larger, if necessary, Lemma 4.3 and transversality guarantee that $F^{\circ n+m}(U_m) \cap (D_{u, m} \times \partial D_s) = \emptyset$. 
 \newline\null\hfilll{\bf QED}
 \enddemo

Set $N = N(d) = m + n$.

\demo {Proof of Theorem 4.1} 
 Now by using the method of the proof of Lemma 4.2 inductively, the same fields of cones

\bigskip\centerline{52}\bigskip\hrule\bigskip
 
$$
 \Cal C_{\bold z} = \bigl\{\,(\xi_1, \xi_2)\bigm| |\xi_2| < |\xi_1|\,\bigr\}
 $$ 
 and 
 $$
 \Cal C_{\bold z}' = \bigl\{\,(\xi_1, \xi_2)\bigm| |\xi_1| < |\xi_2|\,\bigr\}
 $$ 
 are $F^{\circ N}$-trapping and $F^{\circ -N}$-trapping on $U_m$, respectively.
 \newline\null\hfilll{\bf QED}
 \enddemo

\bigskip\centerline{53}\bigskip\hrule\bigskip

\centerline{\bf List Of References}\bigskip

\hangpar\noindent
\hbox to  0.45truein{\bf [A] \hfill}Alekseev, V., {\it Quasirandom dynamical systems I, II, III}, Math. USSR Sbornik, {\bf 5} (1968), pp. 73--128; {\bf 6} (1968), pp. 505--60; {\bf 7} (1969), pp. 1--43.

\hangpar\noindent
\hbox to  0.45truein{\bf [Bi] \hfill}Bieberbach, L., {\it Beispiel zweier ganzer Funktionen zweier komplexer Variabeln, welche eine schlicht volumetreue Abbildung des $\R_4$ auf einen Teil seiner selbst
vermitteln}, Sitzungsber. Preuss. Akad. Wiss. Berlin, Phys.-math. Kl. (1933), pp. 476--79 .

\hangpar\noindent
\hbox to  0.45truein{\bf [Bl] \hfill}Blanchard, P., {\it Complex analytic dynamics on the Riemann sphere}, Bull. (New Ser.) AMS {\bf 11} (1984), pp. 85--141.

\hangpar\noindent
\hbox to  0.45truein{\bf [DN] \hfill}Devaney, R., and Nitecki, Z., {\it Shift automorphisms in the H\'enon mapping}, Comm. math. Phys. {\bf 67} (1979), pp. 137--48.

\hangpar\noindent
\hbox to  0.45truein{\bf [DH] \hfill}Douady, A., and Hubbard, J., {\it On the dynamics of polynomial-like mappings}, Ann. scient. \'Ec. Norm. Sup., $4^e$ ser. {\bf 18} (1985), pp. 287--343.

\hangpar\noindent
\hbox to  0.45truein{\bf [F1] \hfill}Fatou, P., {\it Sur les \'equations fonctionnelles}, Bull. Soc. math. France {\bf 47} (1919), pp. 161--271; {\bf 48} (1920), pp. 33--94; {\bf 48} (1920), pp. 208--314.

\hangpar\noindent
\hbox to  0.45truein{\bf [F2] \hfill}------------, {\it Sur les fonctions m\'eromorphes de deux variables} and {\it Sur certaines fonctions uniformes de deux variables}, C. R. Acad. Sc. Paris {\bf 175}
(1922), pp. 862-65, 1030--33.

\hangpar\noindent
\hbox to  0.45truein{\bf [FM] \hfill}Friedland, S., and Milnor, J., {\it Dynamical properties of plane polynomial automorphisms} (in preparation).

\bigskip\centerline{54}\bigskip\hrule\bigskip

\hangpar\noindent
\hbox to  0.45truein{\bf [H1] \hfill}H\'enon, M., {\it Numerical study of quadratic area preserving mappings}, Q. Appl. Math.  {\bf 27} (1969), pp. 291--312.

\hangpar\noindent
\hbox to  0.45truein{\bf [H2] \hfill}------------, {\it A two-dimensional mapping with a strange attractor}, Commun. math. Phys.  {\bf 50}  (1976), pp. 69--77.

\hangpar\noindent
\hbox to  0.45truein{\bf [HO] \hfill}Hubbard, J., and Oberste-Vorth, R., {\it H\'enon mappings in the complex domain} (in preparation).

\hangpar\noindent
\hbox to  0.45truein{\bf [J] \hfill}Julia, G., {\it M\'emoires sur l'it\'eration des fonctions rationelles}, J. Math. {\bf 8} (1918), pp. 47--245 (See also {\it {\OE}uvres de Gaston Julia},
Gauthier-Villars, Vol. I, Paris, 1968, pp. 121--319.) 

\hangpar\noindent
\hbox to  0.45truein{\bf [M] \hfill}Milnor, J., private communication.

\hangpar\noindent
\hbox to  0.45truein{\bf [Mo] \hfill}Moser, J., {\it Stable and Random Motions in Dynamical Systems}, Princeton Univ. Press, Princeton, N.J., 1973.

\hangpar\noindent
\hbox to  0.45truein{\bf [N] \hfill}Newhouse, S., {\it Lectures on dynamical systems} in {\it Dynamical Systems, C.I.M.E. Lectures, Bressanone, Italy, June 1978}, Birkhauser-Boston, 1980, pp. 1--114.

\hangpar\noindent
\hbox to  0.45truein{\bf [S] \hfill}Smale, S., {\it Diffeomorphisms with many periodic points} (in {\it Differential and Combinatorial Topology}, S. Cairns, ed.), Princeton Univ. Press, Princeton, N.J.,
1965, pp. 63--80. 

\bigskip\centerline{55}\bigskip\hrule\newpage

\centerline{\underbar{\bf ADDENDUM}}
\bigskip
\parindent0pt
I would like to make several remarks:

\parindent10pt\roster
 \item The content of this document is intended to be an unedited copy of the original thesis. A few corrections of minor spelling and other typographical errors have been made. However, this document is
a reconstruction in TeX of the word-processed original. This has been paginated to reduce the number of pages. The original thesis had 65 double spaced pages including the title
page; horizontal lines indicate the original pagination and numbering. 
 \item References to [HO] in the thesis refer to a preprint. In actual publication, the reference on page 2 is to [5], the references on page 5 are generally to [1] and [2], and the reference on page 15 is
to [1].
 \item Crossed mappings---a typical H\'enon-like mapping---are studied also in [2].
 \item The content of this thesis has been published in [3], [4], and [6]:
 \itemitem{$\bullet$} [3] includes the content of chapter 1 and much of chapter 2 using crossed mappings; \newline
 (Note that [3] contains a typographical error: the constant $\frac54 + \frac{\sqrt{5}}2$ of Proposition 2.10 of the thesis is incorrectly given as $\frac54 + \sqrt{\frac52}$ in Proposition 2 of [3].) 
 \itemitem{$\bullet$} [4] includes some of the content of chapter 3; and
 \itemitem{$\bullet$} [6] includes the content of chapters 3 and 4.
 \endroster

\medskip\noindent{\bf References}

\hangpar\noindent
[1] J. Hubbard and R. Oberste-Vorth, {\it H\'enon mappings in the complex domain I: the global topology of dynamical space},	Publ. Math. IHES {\bf 79} (1994), pp. 5-46. 

\hangpar\noindent
[2]	-----------, {\it H\'enon mappings in the complex domain II: projective and inductive limits of polynomials}, in Real \& Complex Dynamical Systems, B. Branner and P. Hjorth, eds., Kluwer, 1995, pp.
89-132.

\hangpar\noindent
[3] R. Oberste-Vorth, {\it Horseshoes among H\'enon mappings}, in {\it Recent Advances in Applied and Theoretical Mathematics}, N. Mastorakis, ed., WSES Press, 2000, pp. 116--121.

\hangpar\noindent
[4] -----------, {\it Horseshoes as projective limits}, in {\it Conference Proceedings: 2002 WSEAS MMACTEE, WAMUS, NOLASC, Vouliagmeni, Athens, Greece,	Dec. 29-31, 2002} (CD-ROM), N. Mastorakis, M. Er, C.
D'Attelis, eds., WSEAS Press, 2003.

\hangpar\noindent
[5] -----------, {\it Normal forms and Fatou-Bieberbach domains}, WSEAS Trans. on Math. {\bf 3} (2004), pp. 253-258.

\hangpar\noindent
[6] -----------, {\it Complex horseshoes}, in {\it Proceedings of the Sixth WSEAS International Conference on Applied Mathematics} (CD-ROM), N. Mastorakis, ed., WSEAS Press, 2004.

\bigskip\null\noindent
Ralph W. Oberste-Vorth\newline
Department of Mathematics\newline
Marshall University\newline
One John Marshall Drive\newline
Huntington, WV 25755\newline
(304) 696-6010 (office)\newline
(304) 696-4646 (fax)\newline
oberstevorth{\@}marshall.edu\newline
\newline
July 1, 2005

\enddocument